\documentclass[letterpaper, 10 pt, conference]{ieeeconf}  
\usepackage{graphicx} 
\usepackage[cmex10]{amsmath} 
\usepackage{amssymb}  
\usepackage[english]{babel}
\usepackage{amsthm}
\usepackage{thmtools}
\usepackage{dsfont} 
\usepackage[usenames,dvipsnames]{xcolor}
\usepackage{hyperref} %
\usepackage{cleveref}
\let\labelindent\relax
\usepackage{enumitem} \setlist[itemize,1]{leftmargin=\dimexpr 25pt}
\usepackage{cite}
\usepackage{url}
\usepackage{microtype}
\usepackage{BOONDOX-ds}
\usepackage{pgfplots} \pgfplotsset{compat=1.14} \usepgfplotslibrary{colorbrewer}
\usepackage{tikz} 
\usepackage{quoting}
\usepackage{siunitx}

\renewcommand{\theenumi}{(\roman{enumi})}

\makeatletter
\def\@IEEEsectpunct{.\ \,}
\def\paragraph{\@startsection{paragraph}{4}{\z@}{1.5ex plus 1.5ex minus 0.5ex}%
{0ex}{\normalfont\normalsize\itshape}}
\makeatother

\declaretheorem[style=definition]{theorem}

\declaretheorem[style=definition]{lemma}
\declaretheorem[style=definition,qed=$\vartriangle$]{remark}

\declaretheorem[style=definition,numbered=no,qed=$\lrcorner$]{standing assumption}
\declaretheorem[style=definition,qed=$\lrcorner$]{assumption}
\declaretheorem[style=definition,qed=$\lrcorner$]{definition}
\declaretheorem[style=definition]{problem}

\newcommand {\nn}{\nonumber}
\newcommand{\beq}{\begin{equation}}
\newcommand{\eeq}{\end{equation}}
\newcommand {\bseq}{\begin{subequations}}
\newcommand {\eseq}{\end{subequations}}
\newcommand {\bma}{\left[}
\newcommand {\ema}{\right]}

\newcommand {\Zplus}{\mathbb{Z}_{+}} 
\newcommand {\R}{\mathbb{R}} 	
 	
\newcommand {\Co}{\mathbb{C}} 	
\newcommand {\Coplus}{\Co_{+}} 	
\newcommand {\Cominus}{\mathbb{C}_{-}} 	
\newcommand {\Cozero}{\mathbb{C}_{0}}

\newcommand{\rank}{\operatorname{rank}} 
\newcommand{\transpose}{\mathsf{T}} 
 
\newcommand{\norm}[1]{\left\lVert#1\right\rVert}

\newcommand{\diag}{\operatorname{diag}}
\newcommand{\spectrum}[1]{{\sigma({#1})}}

\def\eblue{\protect\normalcolor}

\title{On model reduction by least squares moment matching}

\author{Alberto Padoan}

\begin{document}

\maketitle
\thispagestyle{plain}
\pagestyle{plain}

\begin{abstract}     
\noindent 
The paper addresses the model reduction problem by least squares moment matching for continuous-time, linear, time-invariant systems. The basic idea behind least squares moment matching is to approximate a transfer function by ensuring that the interpolation conditions imposed by moment matching are satisfied in a least squares sense.  
This idea is revisited using invariance equations and steady-state responses to provide a new time-domain characterization of least squares moment matching. The characterization, in turn, is then used to obtain a parameterized family of models achieving least squares moment matching. 
The theory is illustrated by a worked-out numerical example. 
\end{abstract}

\section{Introduction}

Model reduction is a central problem in control theory~\cite{antoulas2005approximation}.  Model reduction   occurs frequently in control engineering practice and can be posed,  mathematically, as an approximation problem. For linear time-invariant systems, a popular way to approach the model reduction problem is through methods based on moment matching~\cite{antoulas2005approximation,georgiou1983partial,kimura1986positive,antoulas1990solution,georgiou1999interpolation,grimme1997krylov, gallivan2004sylvester,gallivan2006model}. The main idea is to use (rational) interpolation theory to approximate the transfer function of a system by another transfer function of lower  complexity.   Moment matching consists in imposing that the moments, \textit{i.e.}  the coefficients of the Laurent series expansion, of both transfer functions coincide at given points of the complex plane. Methods based on moment matching are numerically reliable and can be implemented efficiently using Krylov projectors~\cite[Chapter 11]{antoulas2005approximation}. Over the past two decades, moments of linear time-invariant systems have been characterized  in the time domain  using Sylvester equations~\cite{gallivan2004sylvester,gallivan2006model} and, under certain assumptions, using steady-state responses~\cite{astolfi2010model}. These characterizations, in turn, have led to novel model reduction methods by moment matching for nonlinear systems~\cite{astolfi2010model}, for time-delay systems~\cite{scarciotti2016model}, and for systems with isolated singularities~\cite{padoan2017poles,padoan2017mrp,padoan2017eigenvalues,padoan2019isolated1}.  

However, a significant limitation of methods based on moment matching is that the interpolation conditions imposed by moment matching need to hold \emph{exactly}, which, for some purposes, is an unnecessarily stringent assumption. In practice, one can often tolerate a small error around each interpolation point and seek for the ``best'' model which minimizes these errors.  Furthermore, methods based on moment matching typically do not offer  error bounds, which precludes any \textit{a priori} guarantee on the quality of approximation.  

Least squares moment matching provides a particularly interesting solution to both issues~\cite{shoji1985model,aguirre1992least,smith1995least,
aguirre1995algorithm,gugercin2006model,gu2010model}.  The main idea is to require that the interpolation conditions imposed by moment matching are satisfied only in a least squares sense. Model reduction methods based on least squares moment matching  thus overcome the issues mentioned above by minimizing an optimization criterion, which directly yields \textit{a priori} error bounds and, in some cases, guaranteed stability properties~\cite{gugercin2006model,gu2010model}. 
Model reduction by least squares moment matching has a long history~\cite{shoji1985model,aguirre1992least,smith1995least,
aguirre1995algorithm,gugercin2006model,gu2010model,gustavsen1999rational,
mayo2007framework,berljafa2015generalized,berljafa2017rkfit,nakatsukasa2018aaa,antoulas2020interpolatory}, with deep connections with  Pad\'e approximation~\cite{aguirre1992least,aguirre1994model} and Prony's method for filter design~\cite{gugercin2006model,parks1987digital}. The reader is referred to~\cite{antoulas2005approximation,antoulas2020interpolatory}, and references therein, for further detail.

The present work provides a new time-domain characterization of least squares moment matching, which relies on the solution of a constrained optimization problem involving a Sylvester equation.
The main ingredient of our approach is the formalism introduced in~\cite{astolfi2010model}, where the 
moments of a system have been characterized using tools from output regulation theory~\cite{isidori1990output} (see also \cite[Chapter 8]{isidori1995nonlinear}). We show that models achieving least squares moment matching minimize the worst case r.m.s. gain of an error system with respect to the family of output signals produced by a signal generator. Furthermore, we present a new parameterized family of models achieving least squares moment matching. 
Our results offer a unique time-domain perspective on least squares moment matching in terms of invariance equations and steady-state responses, without relying on frequency-domain notions or any other strictly linear tools.  This, in turn, is instrumental to develop a nonlinear enhancement of least squares moment matching, which  will be discussed in a separate publication~\cite{padoan2021lsmr}.

The remainder of this work is organized as follows. 
Section~\ref{sec:problem-formulation-linear} formulates the model reduction problem by least squares moment matching.
Section~\ref{sec:preliminaries} recalls background material on the connection between moment matching
and Sylvester equations~\cite{gallivan2004sylvester,gallivan2006model,astolfi2010model}.
Section~\ref{sec:main-results} presents the main results of the paper.
Section~\ref{sec:examples} illustrates the theory by means of a worked-out example.
Section~\ref{sec:conclusion} concludes the paper with a summary and an outlook to future research directions.   

\textbf{Notation}
$\Zplus$, $\R$ and $\Co$ denote the set of non-negative integers, of real numbers, and of complex numbers, respectively. 
$\Cominus$, $\Cozero$, and $\Coplus$ denote the set of complex numbers with negative real part, zero real part, and positive real part, respectively.
$\iota$ denotes the imaginary unit. 
$e_k$ denotes the vector  with the $k$-th entry equal to one and all other entries equal to zero. 
$I$ denotes the identity matrix.
$J_0$ denotes the matrix with ones on the superdiagonal and zeros elsewhere. 
$J_{s^{\star}}$ denotes the Jordan block associated with the eigenvalue ${s^{\star}\in\Co}$, \emph{i.e.} ${J_{s^{\star}} = s^{\star} I + J_0}$.
$\spectrum{A}$ denotes the spectrum of the matrix ${A \in \R^{n \times n}}$.
$M^{\transpose}$, $M^{\dagger}$ and $\ker M$ denote the transpose, the Moore-Penrose inverse and the kernel of the matrix ${M \in \R^{p \times m}}$, respectively. 
$\norm{\,\cdot\,}_2$ and $\norm{\,\cdot\,}_{2*}$ denote the Euclidean norm on $\R^{n}$ and the corresponding dual norm~\cite[p.637]{boyd2004convex}, respectively. Finally, $f^{(k)}$ denotes the derivative of order ${k\in\Zplus}$ of the function $f$, provided it exists, and $f^{(0)} = f$ by convention. %

\section{Problem formulation} \label{sec:problem-formulation-linear}

Consider a continuous-time, single-input, single-output, linear, time-invariant system described by the equations
\beq \label{eq:system-linear}
\quad \dot{x} = Ax+Bu, \quad y=Cx,
\eeq
in which ${x(t)\in\R^n}$, ${u(t)\in\R}$, ${y(t)\in\R}$ and ${A\in\R^{n \times n}}$, ${B\in\R^{n \times 1}}$ and ${C\in\R^{1\times n}}$ are constant matrices, with transfer function defined as 
$${W(s)=C(sI-A)^{-1}B. }$$
For the notion of moment to make sense, we make the following standing assumption throughout the paper.
\begin{standing assumption}
The system~\eqref{eq:system-linear} is minimal, \textit{i.e.} controllable and observable.
\end{standing assumption}

\begin{definition}   \cite[p.345]{antoulas2005approximation}    \label{def:moment}
The \emph{moment of order ${k \in \Zplus}$} of system~\eqref{eq:system-linear} at ${s^{\star} \in \Co}$, with ${s^{\star} \not \in  \spectrum{A}}$, is defined as the complex number
\beq \nn
\eta_k(s^{\star}) = \frac{(-1)^k}{k!} W^{(k)}(s^{\star}) .
\eeq 
\end{definition}

\noindent 
Given distinct \emph{interpolation points} ${\{s_i\}_{i=1}^N}$, with 
${s_i \in \Co}$ and ${s_i \not \in  \spectrum{A}}$,
and the corresponding \emph{orders of interpolation}
${\{k_i\}_{i=1}^N}$, with ${k_i\in \Zplus}$,
model reduction by moment matching  consists in finding a system
\beq \label{eq:system-rom}
  \quad  \dot{\xi} = F\xi+Gv, \quad \psi = H\xi,
\eeq 
where ${\xi(t) \in\R^{r}}$, ${v(t)\in\R}$, ${\psi(t)\in\R}$ and ${F\in\R^{r \times r}}$, ${G\in\R^{r \times 1}}$ and ${H\in\R^{1\times r}}$ are constant matrices, the transfer function of which
\beq \nn
\hat{W}(s)=H(sI-F)^{-1}G 
\eeq
satisfies the \emph{interpolation conditions}
\beq  \label{eq:matching-condition-linear}
\eta_{j}(s_i)= \hat{\eta}_{j}(s_i) , \quad 0 \le j \le k_i , \quad 1 \le i \le N ,
\eeq 
where $\eta_{j}(s_i) $ and $ \hat{\eta}_{j}(s_i)$ denote the moments of order $j$ of the systems~\eqref{eq:system-linear} and~\eqref{eq:system-rom} at $s_i$, respectively. The system~\eqref{eq:system-rom} is referred to as a \emph{model of system~\eqref{eq:system-linear}} and is said to \emph{achieve moment matching (at $\{ s_i\}^{N}_{i=1}$)} if the interpolation conditions~\eqref{eq:matching-condition-linear} hold~\cite[Chapter 11]{antoulas2005approximation}. Furthermore, if ${r < n}$, then~\eqref{eq:system-rom} is said to be a \emph{reduced order model of system~\eqref{eq:system-linear}}.

\begin{figure}[t]
\centering
\includegraphics[width=\columnwidth,  height=5cm, keepaspectratio]{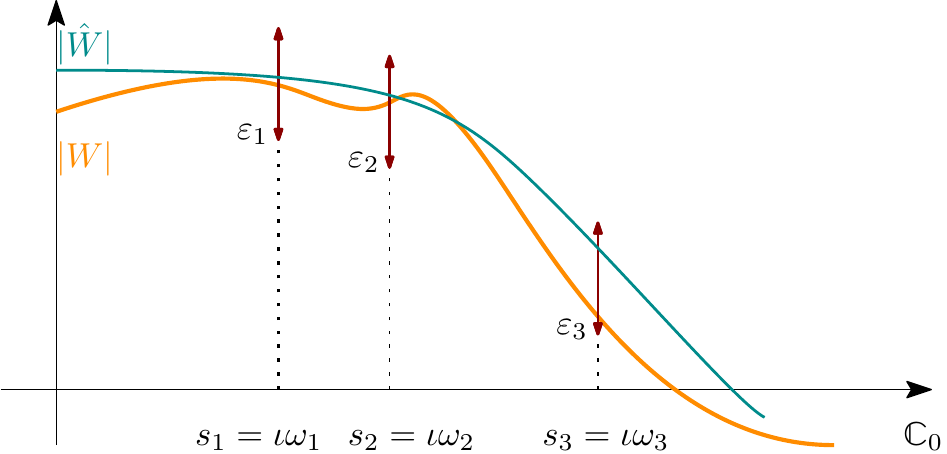}
\caption{The system~\eqref{eq:system-linear} must be approximated around every interpolation point $s_i =\iota\omega_i$ and a small error (depicted as $\varepsilon_i$ for the sake of illustration) can be tolerated. The order $r=1$ of the model is given and the original system is approximated at $N=3$ interpolation points.}
\label{fig:1}
\end{figure}%

We are interested in the following problem. Suppose we wish to approximate system~\eqref{eq:system-linear} with the model~\eqref{eq:system-rom} 
by moment matching at a given set of interpolation points  $\{ s_i\}^{N}_{i=1}$ with corresponding orders of interpolation $\{ k_i\}^{N}_{i=1}$. Assume that the number of interpolation conditions 
$\nu = \sum_{i=1}^{N} (k_i+1)$
is much larger than the order of system~\eqref{eq:system-rom}, \emph{i.e.} $\nu \gg r$, and that  around every interpolation point $s_i$ a small error can be tolerated. The goal is to construct the ``best'' model of system~\eqref{eq:system-linear} of order $r$ (in a sense to be made precise in the sequel).  Fig.~\ref{fig:1}  provides a pictorial representation of the problem, where a model of order $r=1$ needs to be constructed to approximate the original system at $N=3$ interpolation points.  

It is well-known that a model of order $r$ can match up to $2r$ moments~\cite[Chapter 11]{antoulas2005approximation}. The number of interpolation conditions is thus larger than the number of moments that can be matched when $\nu > 2r$. In this case, the interpolation conditions~\eqref{eq:matching-condition-linear} give rise to an overdetermined system of equations which can be solved in a least squares sense. This leads directly to the \emph{model reduction problem by least squares moment matching}.

\begin{problem}
Consider system~\eqref{eq:system-linear}. Let $\{ s_i\}^{N}_{i=1}$ be a set of interpolation points, with ${s_i \in \Co}$ and ${s_i \not \in  \spectrum{A}}$, and let $\{ k_i\}^{N}_{i=1}$ be the corresponding orders of interpolation, with ${k_i\in \Zplus}$. Let ${\nu = \sum_{i=1}^{N} (k_i+1)}$ and ${r\in \Zplus}$, with ${2r < \nu}$.  Find, if possible, a system~\eqref{eq:system-rom} of order $r$ which minimizes the index
\beq \label{eq:index}
\mathcal{J} = \sum_{i=1}^{N} \sum_{j=0}^{k_i} \left| \eta_{j}(s_i) - \hat{\eta}_{j}(s_i) \right|^2 .
\eeq 
\end{problem}

The model~\eqref{eq:system-rom} is said to \emph{achieve least squares moment matching (at $\{ s_i\}^{N}_{i=1}$)} if it minimizes the index~\eqref{eq:index}. Clearly, any model achieving moment matching also achieves least squares moment matching, since the index~\eqref{eq:index} is minimized if the interpolation conditions~\eqref{eq:matching-condition-linear} hold. The set of (reduced order) models of system~\eqref{eq:system-linear} achieving moment matching is therefore a  strict   subset of the set of (reduced order) models of system~\eqref{eq:system-linear} achieving least squares moment matching.

\section{Preliminaries} \label{sec:preliminaries}

We now recall some preliminary results on the connection between moment matching and Sylvester equations from~\cite{gallivan2004sylvester,gallivan2006model,astolfi2010model} with minor variations.

\begin{lemma} \cite[Lemma 2]{astolfi2010model}  \label{lemma:astolfi-2}
Consider system~\eqref{eq:system-linear}.  Let ${k \in \Zplus}$ and ${s^{\star} \in \Co}$, with ${s^{\star} \not \in  \spectrum{A}}$. Then
$${
C\Pi \Psi 
=
\left[ \, \eta_0(s^{\star}) 	\; \eta_1(s^{\star}) \; \cdots \; \eta_k(s^{\star}) \, \right],
}$$
in which ${\Psi\in \R^{(k+1) \times (k+1)}}$ is a signature matrix\footnote{A signature matrix is a diagonal matrix with $\pm 1$ on the main diagonal~\cite[p.44]{bapat2010graphs}.} and ${\Pi   \in \R^{n\times (k+1) }}$ is the unique solution of the Sylvester equation
$${
A\Pi  +B e_1^{\transpose} =\Pi   J_{s^{\star}} .
}$$
\end{lemma}

\begin{lemma} \cite[Lemmas 3 and 4]{astolfi2010model}   \label{lemma:astolfi-3}
Consider system~\eqref{eq:system-linear}. 
Assume ${S\in \R^{\nu \times \nu}}$ is a non-derogatory\footnote{A matrix is non-derogatory if its characteristic  polynomial and its  minimal polynomials coincide \cite[p.178]{horn1994matrix}.} matrix with characteristic polynomial 
\beq \label{eq:characteristic-polynomial}
\chi_S(s)=\prod_{i=1}^{N} (s-s_i)^{k_i+1}
\eeq
and ${L\in\R^{1\times \nu}}$ is such that the pair $(S,L)$ is observable.
Then the moments $\eta_0(s_1),$ $\dots,$ $\eta_{k_1}(s_1)$, $\dots,$ $\eta_0(s_N),$ $\dots,$ $\eta_{k_N}(s_N)$ are in one-to-one correspondence\footnote{The terminology is borrowed from~\cite{astolfi2010model}, where one-to-one correspondence between two objects means that one uniquely determines the other and \emph{vice versa}.} with the matrix ${C\Pi}$, where  ${\Pi \in \R^{n\times \nu}}$ is the unique solution of the Sylvester equation
\beq \label{eq:Sylvester-equation-astolfi}
A\Pi +BL=\Pi S.
\eeq
\end{lemma} %
\noindent
Lemma~\ref{lemma:astolfi-3} establishes that the moments of system~\eqref{eq:system-linear} can be \emph{equivalently} characterized in terms of the product of the output matrix of the system and the solution of the Sylvester equation~\eqref{eq:Sylvester-equation-astolfi}.  This, in turn, is instrumental to provide a time-domain characterization of the moments of system~\eqref{eq:system-linear} in terms of the \emph{steady-state output response}\footnote{See~\cite{isidori1990output} and \cite[Chapter 8]{isidori1995nonlinear}.} of the interconnection of system~\eqref{eq:system-linear} with a signal generator described by the equations
\beq \label{eq:system-signal-generator}
 \dot{\omega} = S\omega, \quad  \theta = L\omega,  
\eeq 
with ${\omega(t) \in \R^\nu}$ and ${\theta(t) \in \R}$, which satisfies the following assumptions.

\begin{assumption} \label{ass:signal-generator-linear}
The matrix ${S\in \R^{\nu \times \nu}}$ is non-derogatory and has characteristic polynomial~\eqref{eq:characteristic-polynomial}. The matrix ${L\in\R^{1\times \nu}}$ is such that the pair $(S,L)$ is observable.
\end{assumption}

\begin{assumption} \label{ass:excitability-linear}
The vector ${\omega(0) \in\R^{\nu}}$ is such that the pair $(S,\omega(0))$ is controllable\footnote{For linear, time-invariant systems, the notions of controllability and excitability of the pair $(S,\omega(0))$ are equivalent. See~\cite{padoan2016geometric} for further detail.}.
\end{assumption}

\begin{theorem} \cite[Theorem 1]{astolfi2010model}   \label{thm:astolfi-1}
Consider system~\eqref{eq:system-linear} and the signal generator~\eqref{eq:system-signal-generator}.
Suppose Assumptions~\ref{ass:signal-generator-linear} and~\ref{ass:excitability-linear} hold.
Assume ${\spectrum{A} \subset \Cominus}$ and ${\spectrum{S} \subset \Cozero}$. Then the moments $\eta_0(s_1),$ $\dots,$ $\eta_{k_1}(s_1)$, $\dots,$ $\eta_0(s_N),$ $\dots,$ $\eta_{k_N}(s_N)$ are in one-to-one correspondence with  the steady-state output response of the interconnected system~\eqref{eq:system-linear}-\eqref{eq:system-signal-generator}, with ${u =  \theta}$.
\end{theorem}

\noindent
Theorem~\ref{thm:astolfi-1} motivates the following notion of model achieving moment matching.

\begin{definition} \cite[p.4]{astolfi2010model} \label{def:reduced-order-model-linear}
The system~\eqref{eq:system-rom} is a \emph{model of system~\eqref{eq:system-linear} at $(S,L)$}, with $S \in \R^{\nu\times \nu}$ such that $\spectrum{S} \cap \spectrum{A} = \emptyset$, if $\spectrum{S} \cap \spectrum{F} = \emptyset$
and 
\beq \label{eq:model-condition-2}
C\Pi = HP,
\eeq
where ${\Pi \in \R^{ n \times \nu}}$ is the unique solution of the Sylvester equation
\eqref{eq:Sylvester-equation-astolfi}, with $L \in \R^{1 \times \nu}$ such that the pair $(S,L)$ is observable, and $P \in \R^{ r \times \nu}$ is the unique solution of the Sylvester equation
\beq \label{eq:Sylvester-equation-astolfi-model} 
FP +GL= P S.
\eeq
In this case, system~\eqref{eq:system-rom} is said to \emph{match the moment of system~\eqref{eq:system-linear} (or to achieve moment matching) at $(S,L)$}. Furthermore, system~\eqref{eq:system-rom} is a \emph{reduced order model of system~\eqref{eq:system-linear} at $(S,L)$} if ${r<n}$. 
\end{definition}

\noindent
A family of reduced order models achieving moment matching for system~\eqref{eq:system-linear} has been determined in~\cite{astolfi2010model} by selecting $r=\nu$ and ${P=I}$ in~\eqref{eq:model-condition-2} and~\eqref{eq:Sylvester-equation-astolfi-model},  respectively.   As a matter of fact, this yields a family of reduced order models of system~\eqref{eq:system-linear} at $(S,L)$ described by equations~\eqref{eq:system-rom}, with
\beq \label{eq:family-linear-astolfi}
F=S-\Delta L, 
\quad G=\Delta, 
\quad  H = C \Pi,
\eeq
in which ${\Delta \in \R^{r \times 1}}$ is such that $\spectrum{S} \cap \spectrum{S-\Delta L} = \emptyset$. The vector ${\Delta }$ is a \emph{parameter} of the family of reduced order models~\eqref{eq:family-linear-astolfi} which can be used to assign prescribed properties to the reduced order model, including stability, passivity, and a given $L_2$-gain~\cite{astolfi2010model}.

\section{Main results} \label{sec:main-results}

We begin our analysis by establishing that models achieving least squares moment matching can be \textit{equivalently} characterized in terms of the solutions of a constrained optimization problem of the form
\beq \label{eq:optimization-problem}
\begin{array}{ll}
    \mbox{minimize}      & \norm{(C \Pi - H P)T }_{2*}^2 \\
    \mbox{subject to}    & FP +GL= P S , \\
    						& \spectrum{S} \cap \spectrum{F} = \emptyset,
\end{array}
\eeq
for  a given  non-singular matrix ${T \in \R^{\nu\times \nu}}$, where ${F\in\R^{r \times r}}$, ${G\in\R^{r \times 1}}$, ${H\in\R^{1\times r}}$ and ${P\in\R^{r\times \nu}}$  are the optimization variables, while system~\eqref{eq:system-linear} and the signal generator~\eqref{eq:system-signal-generator} 
(and, thus, the solution ${\Pi \in \R^{ n \times \nu}}$ of the Sylvester equation~\eqref{eq:Sylvester-equation-astolfi}) 
are problem data. To this end, we first introduce a basic assumption and prove a preliminary lemma, which allows us to rewrite the index~\eqref{eq:index} in terms of the solutions of the Sylvester equations~\eqref{eq:Sylvester-equation-astolfi} and~\eqref{eq:Sylvester-equation-astolfi-model}.
 
\begin{assumption} \label{ass:T-linear}
The matrix ${T \in \R^{\nu\times \nu}}$ is non-singular and such that
\beq \label{eq:thm-index-proof-00}
ST = TJ, \quad L T = \Lambda ,
\eeq
with ${J = \diag(J_{s_1}, \ldots, J_{s_N})}$ and $\Lambda = [\,e_1^{\transpose} \, \cdots \, e_1^{\transpose} \,]$. 
\end{assumption}

\begin{remark} \label{rem:T-unitary}
Assumption~\ref{ass:signal-generator-linear} implies that the pair $(S,L)$ is observable and, thus, it guarantees the existence of a matrix ${T \in \R^{\nu\times \nu}}$ such that Assumption~\ref{ass:T-linear} holds. 
Note that, without loss of generality, the matrix $T$ can be assumed to be unitary\footnote{A matrix ${U \in \R^{n \times n}}$ is unitary if ${UU^{\transpose} =  I}$
\cite[p.84]{horn1994matrix}.} if the matrix $S$ is normal\footnote{A real matrix ${N \in \R^{n \times n}}$ is normal if it commutes with its transpose, \textit{i.e.} ${AA^{\transpose} =  A^{\transpose}A}$~\cite[p.131]{horn1994matrix}.} and $\norm{L}_{2*}=1$.
\end{remark}

\begin{lemma} \label{lem:index}
Consider system~\eqref{eq:system-linear}, the model~\eqref{eq:system-rom} and the signal generator~\eqref{eq:system-signal-generator}.  Suppose Assumptions~\ref{ass:signal-generator-linear} and~\ref{ass:T-linear} hold. Assume ${\spectrum{S} \cap \spectrum{A} = \emptyset}$ and ${\spectrum{S} \cap \spectrum{F} = \emptyset}$. Then 
\beq \label{eq:thm-index-inequality}
\norm{(C \Pi - H P)T }_{2*}^2 
=
\sum_{i=1}^{N} \sum_{j=0}^{k_i} \left| \eta_{j}(s_i) - \hat{\eta}_{j}(s_i) \right|^2 ,
\eeq
where ${\Pi \in\R^{n \times \nu}}$ and ${P\in\R^{r \times \nu}}$ are the (unique) solutions of the Sylvester equations~\eqref{eq:Sylvester-equation-astolfi} and~\eqref{eq:Sylvester-equation-astolfi-model}, respectively.
\end{lemma}

\begin{proof}  
The assumptions $\spectrum{S} \cap \spectrum{A} = \emptyset$ and $\spectrum{S} \cap \spectrum{F} = \emptyset$ directly imply existence and uniqueness of the solutions of the Sylvester equations~\eqref{eq:Sylvester-equation-astolfi} and~\eqref{eq:Sylvester-equation-astolfi-model}, respectively. Furthermore, by Assumptions~\ref{ass:signal-generator-linear} and~\ref{ass:T-linear}, \eqref{eq:Sylvester-equation-astolfi} and~\eqref{eq:thm-index-proof-00} together imply
\beq \nn
A\Pi +B\Lambda T^{-1}=\Pi TJT^{-1},
\eeq
or, equivalently,
\beq \nn
A\Pi T +B\Lambda =\Pi TJ ,
\eeq
where ${\Pi \in \R^{n \times \nu}}$ is the (unique) solution of the Sylvester equation~\eqref{eq:Sylvester-equation-astolfi}. Then, setting  ${\widetilde \Pi = \Pi T}$ and appealing to Lemma~\ref{lemma:astolfi-2}, one obtains
\beq \label{eq:thm-index-proof-01}
C\widetilde \Pi \Psi \!=\! \left[ \, \eta_0(s_1) 	\,  \cdots \, \eta_{k_1}(s_1) \, \cdots \, \eta_0(s_N) 	\,  \cdots \, \eta_{k_N}(s_N) \, \right]  \!,
\eeq
in which ${\Psi \in \R^{\nu \times \nu}}$ is a signature matrix. A similar reasoning applies to the model~\eqref{eq:system-rom}. Thus, setting ${\widetilde P = P T}$ and appealing to Lemma~\ref{lemma:astolfi-2}, yields 
\beq \label{eq:thm-index-proof-02}
H\widetilde P  \Psi \!=\! \left[ \, \hat\eta_0(s_1) 	\, \cdots \, \hat\eta_{k_1}(s_1) \, \cdots \, \hat\eta_0(s_N) 	\,  \cdots \, \hat\eta_{k_N}(s_N) \, \right] \!,
\eeq 
where ${P\in \R^{r \times \nu}}$ is the (unique) solution of the Sylvester equation~\eqref{eq:Sylvester-equation-astolfi-model}.
Then 
\beq \nn
\scalebox{0.95}{$
\begin{array}{rcl}
\norm{(C \Pi - H P)T }_{2*}^2 
&=& \norm{C\widetilde \Pi  - H \widetilde P }_{2*}^2  \\
&=& \norm{C\widetilde \Pi \Psi - H \widetilde P \Psi }_{2*}^2  \\
&\stackrel{\eqref{eq:thm-index-proof-01},\eqref{eq:thm-index-proof-02}}{=}& 
\displaystyle
\sum_{i=1}^{N} \sum_{j=0}^{k_i} \left| \eta_{j}(s_i) - \hat{\eta}_{j}(s_i) \right|^2 ,
\end{array}
$}
\eeq
where the second identity holds since $\Psi$ is a signature matrix and since the norm $\norm{\,\cdot\,}_{2*}$ is unitarily invariant.
\end{proof}%

We are now ready to show that, under certain assumptions, the interpolation constraints imposed by least squares moment matching can be \textit{equivalently} characterized in terms of the solutions of the optimization problem~\eqref{eq:optimization-problem}.

\begin{theorem} \label{thm:moments-optimization-linear}
Consider system~\eqref{eq:system-linear}, the model~\eqref{eq:system-rom} and the signal generator~\eqref{eq:system-signal-generator}.  Suppose Assumptions~\ref{ass:signal-generator-linear} and~\ref{ass:T-linear} hold. Assume ${\spectrum{S} \cap \spectrum{A} = \emptyset}$. Then the model~\eqref{eq:system-rom} achieves least squares moment matching at $\spectrum{S}$ if and only if there exists a full rank matrix ${P\in\R^{r\times \nu}}$ such that $(F,G,H,P)$ is a solution of the optimization problem~\eqref{eq:optimization-problem}.
\end{theorem}

\begin{proof}
``$\Rightarrow$''. Suppose the model~\eqref{eq:system-rom} achieves least squares moment matching at $\spectrum{S}$. By definition, the moments of the model~\eqref{eq:system-rom} at ${s_i}$ up to the order ${k_i}$ are well-defined for every ${i  \in \{1, \ldots,  N  \}}$. Then ${\spectrum{S} \cap \spectrum{F} = \emptyset}$, which implies existence and uniqueness of a matrix $P$ which solves the Sylvester equation~\eqref{eq:Sylvester-equation-astolfi-model}.  We now show that the matrix $P$ is full rank. To this end, assume, by contradiction, that $(F,G)$ is not reachable. Then the Popov-Belevitch-Hautus criterion implies that there exists ${\lambda \in \Co}$ and ${w\in\Co^r}$, with ${w\not  = 0}$, such that $w^{\transpose} [\, \lambda I - F \, | \, G \,] = 0$. Equivalently, $\lambda w^{\transpose} = w^{\transpose}  F$ and $w^{\transpose}G = 0$.  Then 
\beq \nn
w^{\transpose} P S 
\stackrel{\eqref{eq:Sylvester-equation-astolfi-model} }{=}
   w^{\transpose} (FP+GL) 
 =  w^{\transpose}FP+w^{\transpose}GL 
 =  \lambda w^{\transpose}P,
\eeq
which yields the contradiction $\spectrum{S} \cap \spectrum{F} \not = \emptyset$. Then $\rank P = r$ by \cite[Theorem 1]{de1981controllability}, since $(F,G)$ is reachable and $(S,L)$ is observable  by Assumption~\ref{ass:signal-generator-linear}.  We conclude that $(F,G,H,P)$ satisfies the constraints of~\eqref{eq:optimization-problem}. Furthermore, $(F,G,H,P)$ minimizes the objective function of~\eqref{eq:optimization-problem} by Lemma~\ref{lem:index} and by the assumption that~\eqref{eq:system-rom} achieves least squares moment matching at $\spectrum{S}$.

``$\Leftarrow$''. The implication follows directly from Lemma~\ref{lem:index} and from the definition of model achieving least squares moment matching. 
\end{proof}

\begin{remark}
Theorem~\ref{thm:moments-optimization-linear} provides a new characterization of least squares moment matching in terms of the solutions of the optimization problem~\eqref{eq:optimization-problem}. In analogy with the discussion presented in~\cite{astolfi2010model}, this characterization allows one to \textit{define} a time-domain notion of least squares moment matching for systems which do not possess a representation in terms of a transfer function, including nonlinear systems as well as linear time-varying systems. An in-depth analysis of this point is beyond the scope of this paper and will be presented in a separate publication~\cite{padoan2021lsmr}.
\end{remark}

\subsection{An \textit{a priori} error bound}

Least squares moment matching can be given a simple interpretation in terms of the steady-state behavior of the error system 
\beq \label{eq:system-error-linear}
 \dot{x} = Ax+B u,  \quad 
 \dot{\xi} = F\xi+G u,   \quad 
 e=Cx-H\xi,  
\eeq 
in which  ${x(t)\in\R^n}$, ${\xi(t)\in\R^r}$, ${u(t)\in\R}$, and ${e(t)\in\R}$. In particular, if all solutions of the signal generator~\eqref{eq:system-signal-generator} are periodic and if the steady-state output response $e_{ss}$ of the  interconnected system~\eqref{eq:system-signal-generator}-\eqref{eq:system-error-linear}, with $u=\theta$, is well-defined, then achieving least squares moment matching corresponds to minimizing an upper bound of the \textit{worst case r.m.s. gain} of the error system~\eqref{eq:system-error-linear}  with respect  to the family of output signals produced by the signal generator~\eqref{eq:system-signal-generator}, defined as~\cite[p.98]{boyd1991linear}
\beq \label{eq:gain-rms-linear}
\gamma_{rms} =  \sup_{u \in \mathcal{U}} \frac{\norm{e_{ss}}_{rms}}{\norm{u}_{rms}}
\eeq 
where $\norm{v}_{rms}$ is the \textit{r.m.s. value} of the signal ${v(t) \in \R^q}$, defined as~\cite[p.86]{boyd1991linear}
\beq \label{eq:rms}
\norm{v}_{rms} = \left(\lim_{\tau \to \infty} \frac{1}{\tau} \int_{0}^{\tau} v(t)^{\transpose}v(t)dt \right)^{1/2} ,
\eeq
provided the limit exists, while $\mathcal{\mathcal{U}}$ is the family of output signals produced by~\eqref{eq:system-signal-generator} with non-zero r.m.s. value.

\begin{theorem} \label{thm:rms-linear}
Consider system~\eqref{eq:system-linear}, the model~\eqref{eq:system-rom} and the signal generator~\eqref{eq:system-signal-generator}. Suppose Assumptions~\ref{ass:signal-generator-linear},~\ref{ass:excitability-linear} and~\ref{ass:T-linear} hold. Assume ${\spectrum{A} \cup \spectrum{F} \subset \Cominus}$, ${S + S^{\transpose} = 0}$ and $\norm{L}_{2*}=1$. Then the following statements hold.
\begin{itemize}
\item[(i)] The steady-state output response of the interconnected system~\eqref{eq:system-signal-generator}-\eqref{eq:system-error-linear}, with ${u=\theta}$, is well-defined and uniquely determined by the moments of the error system~\eqref{eq:system-error-linear}  at $\spectrum{S}$.
\item[(ii)] The worst case r.m.s. gain of the error system~\eqref{eq:system-error-linear} 
 with respect  to the family of output signals produced by the signal generator~\eqref{eq:system-signal-generator}  is well-defined and such that
\beq \label{eq:gain-bound-linear}
\gamma_{rms}  \le \norm{C \Pi - H P}_{2*} ,
\eeq
with ${\Pi \in\R^{n \times \nu}}$ and ${P\in\R^{r \times \nu}}$ the unique solutions of the Sylvester equations~\eqref{eq:Sylvester-equation-astolfi} and~\eqref{eq:Sylvester-equation-astolfi-model}, respectively. 
\item[(iii)]  The error bound~\eqref{eq:gain-bound-linear} is minimized if the model~\eqref{eq:system-rom} achieves least squares moment matching at~$\spectrum{S}$.
\end{itemize}
\end{theorem}

\begin{proof}
(i). The dynamics of the interconnected system~\eqref{eq:system-signal-generator}-\eqref{eq:system-error-linear}, with $u=\theta$, is governed by the equations
\beq \label{eq:system-error-interconnected-linear} 
\dot{\omega} \!=\! S\omega, \
\dot{x} \!=\! Ax+B L\omega,  \
 \dot{\xi} \!=\! F\xi+G L\omega,   \
 e\!=\!Cx-H\xi .
\eeq 
By assumption, ${\spectrum{A} \cup \spectrum{F} \subset \Cominus}$. Moreover, ${\spectrum{S} \subset \Cozero}$, 
since ${S + S^{\transpose} = 0}$. Then, by the center manifold theorem~\cite[p.4]{carr1981application}, the interconnected system~\eqref{eq:system-error-interconnected-linear} has a well-defined center manifold
\beq \nn
\mathcal{M}_c = \left\{ \, (x,\xi,\omega) \,:\,  x = \Pi\omega, \, \xi =  P\omega \,   \right\} ,
\eeq
where ${\Pi \in\R^{n \times \nu}}$ and ${P\in\R^{r \times \nu}}$ are the (unique) solutions of the Sylvester equations~\eqref{eq:Sylvester-equation-astolfi} and~\eqref{eq:Sylvester-equation-astolfi-model}, respectively. The manifold $\mathcal{M}_c$ is exponentially attractive, since
\beq\nn
\scalebox{0.95}{$
\renewcommand{\arraystretch}{1}
\setlength{\arraycolsep}{1.5 pt}
\dot{\overbrace{\bma \begin{array}{c} x- \Pi\omega \\ \xi - P \omega \end{array} \ema}} 
\!\stackrel{\eqref{eq:system-error-interconnected-linear}}{=} \!
\bma \begin{array}{c} Ax + BL\omega - \Pi S\omega \\ F\xi + GL\omega - P S\omega \end{array} \ema 
\!\stackrel{\eqref{eq:Sylvester-equation-astolfi},\eqref{eq:Sylvester-equation-astolfi-model}}{=}\!
\bma \begin{array}{c} A(x- \Pi\omega) \\ F(\xi - P \omega) \end{array} \ema \!.$}
\eeq
The output of the interconnected system~\eqref{eq:system-error-interconnected-linear} can be thus written as
$e(t) =  e_{ss}(t) + e_{d}(t),$
in which
\beq \label{eq:error-steady-state-linear}
e_{ss}(t) = (C\Pi - HP)\omega(t)
\eeq
is the well-defined steady-state output response of the error system~\eqref{eq:system-error-linear} 
and $e_{d}$ is an exponentially decaying signal. Furthermore, under the stated assumptions, Theorem~\ref{thm:astolfi-1} implies that the steady-state output response $e_{ss}$ is uniquely determined by the moments of the error system~\eqref{eq:system-error-linear} at $\spectrum{S}$.

(ii). The worst case r.m.s. gain of the error system~\eqref{eq:system-error-linear} is  well-defined: the matrix $S$ is skew-symmetric, which implies that the solutions of the signal generator~\eqref{eq:system-signal-generator} are periodic, and the steady-state output response of the interconnected system~\eqref{eq:system-error-interconnected-linear} is also periodic (with the same period) by~\eqref{eq:error-steady-state-linear}.   To obtain the error bound~\eqref{eq:gain-bound-linear}, note that 
\begin{align}
\norm{e_{ss}(t)}_{2}  
&\stackrel{\eqref{eq:error-steady-state-linear}}{=} \!
\norm{(C\Pi - HP)\omega(t)}_{2} \nn \\
&\,\le
\norm{(C\Pi - HP)T}_{2*}
\norm{T^{-1} \omega(t)}_{2} \nn \\
&\,\le
\norm{C\Pi - HP}_{2*}
\norm{L\omega(t)}_{2} , \label{eq:coordinates-change-linear1}
\end{align}
where the second inequality follows for any  non-singular matrix ${T\in\R^{\nu \times \nu}}$
from the Cauchy-Schwartz inequality, while the third inequality follows taking $T$ unitary and such that the first row of $T^{-1}$ is $L$  (which is always possible in view of Assumptions~\ref{ass:signal-generator-linear},~\ref{ass:excitability-linear},~\ref{ass:T-linear} and Remark~\ref{rem:T-unitary}).  Then
\beq \nn
\gamma_{rms} 
\!\stackrel{\eqref{eq:gain-rms-linear}}{=} \! 
\sup_{u \in \mathcal{U}} \frac{\norm{e_{ss}}_{rms}}{\norm{u}_{rms}}
\!\stackrel{\eqref{eq:coordinates-change-linear1}}{\le} 
\norm{C\Pi - HP}_{2*} ,
\eeq
which proves the error bound~\eqref{eq:gain-bound-linear}.

(iii) By Theorem~\ref{thm:moments-optimization-linear}, if the model~\eqref{eq:system-rom} achieves least squares moment matching at $\spectrum{S}$, then $(F,G,H,P)$ is a solution of the optimization problem~\eqref{eq:optimization-problem} for some full rank matrix ${P\in\R^{r\times \nu}}$ and some matrix ${T \in \R^{\nu\times \nu}}$ such that Assumption~\ref{ass:T-linear} holds. Now recall that, by assumption, the matrix $S$ is skew-symmetric and $\norm{L}_{2*}=1$. By Remark~\ref{rem:T-unitary}, the matrix $T$ can be thus taken to be unitary without loss of generality. We conclude that every model~\eqref{eq:system-rom} achieving least squares moment matching at $\spectrum{S}$ minimizes the error bound~\eqref{eq:gain-bound-linear}, since the right-hand-side of~\eqref{eq:gain-bound-linear} coincides with the objective function of the optimization problem~\eqref{eq:optimization-problem}, which for $T$ unitary  becomes $\norm{C\Pi - HP}_{2*}$.
\end{proof}

\subsection{A family of models achieving least squares moment matching}  \label{ssec:LS-MR-MM-linear}

The results of the previous sections can be taken as a starting point to obtain a family of models achieving least squares moment matching by regarding the optimization variable $P$ as a given \emph{parameter}.
The family of models in question is described by the equations~\eqref{eq:system-rom}, with 
\beq \label{eq:family-linear}
F=P(S-\Delta L)Q,
\quad G=P \Delta , 
\quad  H = C \Pi Q,
\eeq
in which ${S\in \R^{\nu \times \nu}}$ is a non-derogatory matrix with characteristic polynomial~\eqref{eq:characteristic-polynomial} such that $\spectrum{S} \cap \spectrum{A} = \emptyset$, 
${L\in\R^{1\times \nu}}$ is such that the pair $(S,L)$ is observable,
${\Pi \in \R^{n \times \nu}}$ is the (unique) solution of the Sylvester equation~\eqref{eq:Sylvester-equation-astolfi},
while ${P \in \R^{r \times \nu}}$, ${\Delta \in\R^{\nu \times 1}}$ and ${Q \in\R^{\nu \times r}}$  are such that
\begin{itemize}
\item[(A$_P$)]  the matrix $P$ is full rank and such that the subspace ${\ker P}$ is  a \eblue  ${(S,L)}$ conditioned invariant,
\item[(A$_Q$)]    the matrix $Q$ is a weighted Moore-Penrose inverse of $P$~\cite[p.117]{ben2003generalized}, defined as  
$${Q =   P^{\dagger}_M = M P^{\transpose} (PM P^{\transpose})^{-1}},   $$
with ${M = TT^{\transpose}}$ and $T$ such that Assumption~\ref{ass:T-linear} holds,
\item[(A$_\Delta$)] the vector $\Delta$ is such that the subspace~$\ker P$ is ${(S-\Delta L)}$-invariant and such that 
$${\spectrum{S} \cap \spectrum{P(S-\Delta L)Q} = \emptyset},$$ 
\end{itemize}
in which case $P$, $\Delta$ and $Q$ are said to be \emph{admissible} for the parameterization~\eqref{eq:family-linear}.

The family of models~\eqref{eq:family-linear} provides a solution to the model reduction problem by least squares moment matching, as detailed by the following statement.    

\begin{theorem} \label{thm:1}
Consider system \eqref{eq:system-linear} and the family of models \eqref{eq:family-linear}. Let ${P\in\R^{r\times \nu}}$,  ${Q\in\R^{\nu\times r}}$ and ${\Delta\in\R^{\nu \times 1}}$ be admissible for the parameterization \eqref{eq:family-linear}. Then the  family of models~\eqref{eq:family-linear}   achieves least squares moment matching at $\spectrum{S}$.
\end{theorem}

\begin{proof}
By Theorem~\ref{thm:moments-optimization-linear}, it suffices to establish that $(F,G,H,P)$, with $F$, $G$ and $H$ defined as in~\eqref{eq:family-linear}, is a solution of the optimization problem~\eqref{eq:optimization-problem} for some full rank matrix ${P \in \R^{r\times \nu}}$ and some  matrix 
 ${T \in \R^{\nu\times \nu}}$ such that Assumption~\ref{ass:T-linear} holds. 
Equivalently,  it suffices to show that all admissible parameters solve the optimization problem
\beq \label{eq:optimization-problem-proof}
\begin{array}{ll}
    \mbox{minimize}      & \norm{C \Pi(I-QP)T}_{2*}^2 \\
    \mbox{subject to}    & P(S-\Delta L)(I-QP) = 0, \\
    						& \spectrum{S} \cap \spectrum{P(S-\Delta L)Q} = \emptyset,
\end{array}
\eeq
with optimization variables $P$, $Q$ and $\Delta$,  since~\eqref{eq:optimization-problem-proof} can be obtained, after simple algebraic manipulations, by direct substitution of~\eqref{eq:family-linear} into~\eqref{eq:optimization-problem}.

Let ${P}$, ${\Delta}$ and ${Q}$ be admissible for the parameterization \eqref{eq:family-linear} and let ${T}$ be such that Assumption~\ref{ass:T-linear} holds. Note that (A$_\Delta$) trivially implies the condition ${\spectrum{S} \cap \spectrum{P(S-\Delta L)Q} = \emptyset}$. Furthermore, (A$_P$) and (A$_\Delta$) imply that the subspace ${\ker P}$ is $(S-\Delta L)$-invariant, \textit{i.e.} 
\beq \label{eq:lemma1-linear-1}
P(S-\Delta L)v =0,   \quad \forall \,  v \in \ker P.  
\eeq
Let ${\bar{v}_i  = (I - QP)e_i}$, with ${i \in \{1, \ldots, \nu\}}$. Note that ${\bar{v}_i \in \ker P}$, since 
${P\bar{v}_i = P(I - QP)e_i = 0}$ in view of ($\text{A}_{Q}$). Then~\eqref{eq:lemma1-linear-1} implies
$P(S-\Delta L)\bar{v}_i = 0,$ with $i  \in \{1, \ldots, \nu\}$, 
or, equivalently,
$
P(S-\Delta L) (I - QP)  =0. 
$
Finally,   recall that $Q$ is a weighted Moore-Penrose inverse of $P$ and, hence, $Q$ is a minimum-norm least-squares solution of the equation ${C\Pi (I - QP) = 0}$~\cite[p.117]{ben2003generalized}.  
\end{proof}

\begin{remark} \label{rem:properties}
The parameters $P$ and $\Delta$ can be used to assign prescribed properties to the reduced order model~\eqref{eq:family-linear}. For example, in order to preserve the first ${r}$ dominant eigenvalues of the original system
(and, hence, its stability properties) the parameters $P$ and $\Delta$ can be selected as follows. The vector $\Delta$
is selected in such a way that the matrix $(S-\Delta L)$ preserves the first ${\nu}$ dominant eigenvalues of the original system. This is always possible since the pair $(S,L)$ is observable. The matrix $P$ is  defined as $P = [\, P_1^{\transpose} \ \cdots \ P_{r}^{\transpose} \,]^{\transpose} ,$
with $\{P_1,\ldots, P_r\}$ a real Jordan basis of right eigenvectors of the matrix $S-\Delta L$ corresponding to the first ${r}$ dominant eigenvalues of the matrix ${(S-\Delta L)}$. By construction, this yields 
\beq \label{eq:properties}
P (S-\Delta L) = F P
\eeq
for some matrix ${F \in \R^{ r \times r}}$. This ensures that the parameters $P$ and $\Delta$ are admissible, since $P$ is full rank and such that the subspace $\ker P$ is both a $(S,L)$ conditioned invariant and ${(S-\Delta L)}$-invariant. Furthermore, \eqref{eq:properties} ensures that the reduced order model preserves the first ${r}$ dominant eigenvalues of the matrix ${(S-\Delta L)}$ and, hence, those of the original system. We illustrate this construction by a worked-out example in Section~\ref{sec:examples}.
\end{remark}

\section{Example}  \label{sec:examples}

We consider the flexible space structure benchmark model from~\cite{gawronski1991model} (see also~\cite{gawronski1996balanced}). The system is described by the equations~\eqref{eq:system-linear}, with
\bseq \label{eq:system-FSS-matrices}
\begin{align}
A &= \diag(A_1,\ldots,A_{K}), \\
B &= [\, B_1^{\transpose} \ \cdots \ B_{K}^{\transpose} \,]^{\transpose}, \\
C &= [\, C_1 \ \cdots \ C_{K} \,],
\end{align}
\eseq
in which the integer ${K>0}$ is the number of modes of the structure and 
\beq \nn
A_k = 
\bma
\begin{array}{cc}
-2\chi_k \phi_k & -\phi_k \\
\phi_k & 0
\end{array}
\ema, \,
B_k = 
\bma
\begin{array}{c}
b_k  \\
0
\end{array}
\ema, \,
C_k = 
\bma
\begin{array}{c}
c_{rk}\\
\tfrac{c_{dk}}{\phi_k}
\end{array}
\ema^{\transpose},
\eeq
where ${\chi_k \in (0,0.001)}$, ${\phi_k \in (0,100)}$,  ${b_k \in (0,1)}$ and ${C_k \in (0,10)^{1\times 2}}$ are uniformly distributed random numbers (generated in MATLAB with the function \texttt{rand} and seed $1009$). The number of modes selected for the simulations is ${K=30}$ and, thus, the order of the   original  system is ${n=60}$.

Suppose we wish to construct a reduced order model that approximates well the original system at low frequency, say below $20$ rad/s. Further, suppose we wish the reduced order model to preserve the dominant eigenvalues of the original system (and, hence, its stability properties).

For illustration, the parameterization~\eqref{eq:family-linear} is used to build a reduced order model of order ${r=10}$  which meets the desired specifications   while achieving   least squares moment matching at $\{\pm \iota\omega_i\}^{12}_{i=0}$, with ${\omega_1 = 0.01}$, ${\omega_2 = 0.1}$, ${\omega_3 = 1}$, ${\omega_4 = 5.5}$, ${\omega_5 = 10}$, ${\omega_6 = 16}$, ${\omega_7 = 20}$, ${\omega_8 = 30}$, ${\omega_9 = 50}$, ${\omega_{10} = 100}$, ${\omega_{11} = 1000}$, and ${\omega_{12} = 10000}$.  To this end, the   matrices $S$ and $L$ are defined as 
\beq \nn
S = \diag(S_1,S_2,\ldots,S_{11},S_{12}), \quad 
L = \tfrac{1}{\sqrt{24}}[ \ \underbrace{1 \  \cdots \   1}_{24} \ ], 
\eeq
where
\beq \nn
S_i = 
\bma
\begin{array}{cc}
0 & \omega_i\\
-\omega_i & 0  
\end{array}
\ema, \hfill \quad 1 \le i \le 12 .
\eeq 
This defines an observable signal generator of order ${\nu=24}$ described by the equations~\eqref{eq:system-signal-generator} such that  $\spectrum{S} = \{\pm \iota\omega_i\}^{12}_{i=0}$.  The Sylvester equation~\eqref{eq:Sylvester-equation-astolfi} is solved (in MATLAB with the function \texttt{sylv}) and the solution $\Pi$ is used to define the matrix ${H = C\Pi}$. A standard pole placement algorithm (implemented in MATLAB by the function \texttt{place}) is used to select the vector $\Delta$ in such a way that $S-\Delta L$ preserves the first ${\nu = 24}$ dominant eigenvalues of the original system. 
The matrix $P$ is defined as
$P = [\, P_1^{\transpose} \ \cdots \ P_{r}^{\transpose} \,]^{\transpose}$,
with $\{P_1,\ldots, P_r\}$ a real Jordan basis of right eigenvectors corresponding to the first ${r = 10}$ dominant eigenvalues of the matrix $S-\Delta L$. The matrix $Q$ is defined as ${Q = P^{\dagger}}$. By Remark~\ref{rem:properties}, this ensures that the parameters $\Delta$, $P$, and $Q$ are admissible and, hence, that the reduced order model model achieves least squares moment matching at $\{\pm \iota\omega_i\}^{12}_{i=0}$, in agreement with Theorem~\ref{thm:moments-optimization-linear}. Furthermore, this ensures that the reduced order model preserves the first ${r = 10}$ dominant eigenvalues of the original system.

\begin{figure}[t!]
\centering
\definecolor{col1}{RGB}{217,95,2}
\definecolor{col2}{RGB}{117,112,179}
\definecolor{col3}{RGB}{27,158,119} 
\pgfplotsset{compat=1.14}
\begin{tikzpicture}
\matrix{
\begin{loglogaxis}[
height = 0.25\textwidth,
width = 0.475\textwidth,
xmin=10^-2,xmax=10^4,
ymin=10^-3,ymax=10^3,
xmajorticks=false,
ylabel={{$ \scriptstyle  |W(\iota\omega)|, \ |\hat{W}(\iota\omega)|$}},
grid
]
\addplot [col3, very thick, smooth]  table [x index = {0}, y index = {1}, col sep=comma]{FSS_sys.csv};
\addplot [col1, very thick, smooth, dashed] table [x index = {0}, y index = {1}, col sep=comma]{FSS_rom.csv};
\end{loglogaxis}  
\\
\begin{loglogaxis}[
height = 0.25\textwidth,
width = 0.475\textwidth,
xmin=10^-2,xmax=10^4,
ymin=10^-3,ymax=10^3,
xlabel={{$ \omega ~\text{[rad/s]} $}},
ylabel={{$ \scriptstyle  |W(\iota\omega) - \hat{W}(\iota\omega)|/|W(\iota\omega)|$}},
grid
]
\addplot [col1, very thick, smooth]  table [x index = {0}, y index = {1}, col sep=comma]{FSS_err_rom.csv};
\end{loglogaxis} 
\\
};
\end{tikzpicture}
\centering
\vspace{-0.3cm}
\caption{Top: Frequency response of  the original system  (solid) and of the reduced order model (dashed). Bottom: Frequency response of the corresponding relative error.
} 
\label{fig:FSS_rom} 
\end{figure}
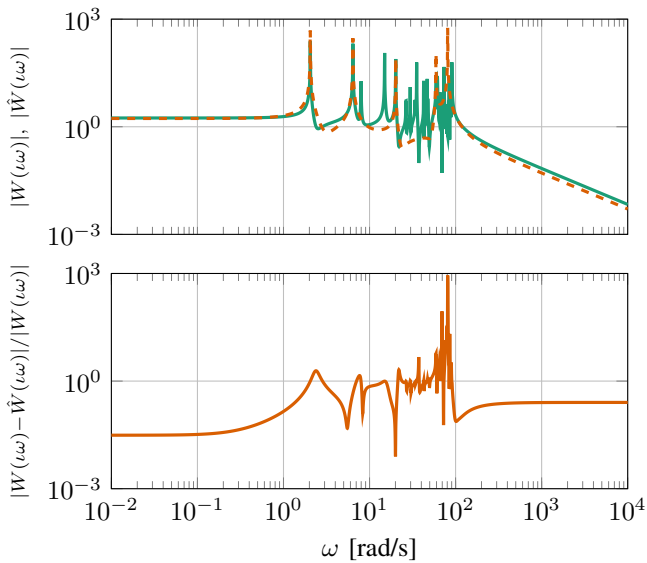

Simulations have been run using standard routines of MATLAB and the experiments have been performed in double precision on a 3.5 GHz Intel Core i7 processor. The top part of Fig.~\ref{fig:FSS_rom} shows the frequency response of the original system (solid) and of the reduced order model (dashed), respectively. The bottom part of Fig.~\ref{fig:FSS_rom} shows the frequency response of the corresponding relative error. 
Note that the reduced order model approximates well the original system at low frequency (below $20$ rad/s) and preserves the first $r$ dominant eigenvalues of the original system, as required. 
Furthermore,  selecting for illustration the initial condition  ${\omega(0)=L^{\transpose}}$ yields ${\norm{e_{ss}}_{r.m.s.} \approx 0.1218}$  (which can be computed in MATLAB with the function \texttt{rms}) and ${\norm{C\Pi - HP}_{2*} \approx  0.5871}$, in agreement with Theorem~\ref{thm:rms-linear}.

\section{Conclusion} \label{sec:conclusion}

The model reduction problem by least squares moment matching has been studied. A new characterization of the notion of least squares moment matching has been presented exploiting invariance equations and steady-state responses. The theory developed does not rely on frequency-domain notions or any other strictly linear tools, thus offering a unique time-domain perspective on least squares moment matching.  The nonlinear enhancement of our results is the subject of ongoing research and will be discussed in a separate publication~\cite{padoan2021lsmr}.

\section*{Acknowledgment}                               
The author warmly thanks Dr. F. Forni for suggesting the problem and Dr. A. Astolfi for his constant support.

\bibliographystyle{IEEEtran}
\bibliography{refs}

\end{document}